\documentclass[twocolumn,aps,showpacs,amssymb,a4paper]{revtex4}
\pdfoutput=1

\usepackage{graphicx,amssymb,amsmath,relsize}
\usepackage{url,color}
\usepackage[pdftex,colorlinks]{hyperref}
\definecolor{darkblue}{cmyk}{1,0,0,0.8}
\definecolor{darkred}{cmyk}{0,1,0,0.7}
\hypersetup{anchorcolor=black,
  citecolor=darkblue, filecolor=darkblue,
  menucolor=darkblue,pagecolor=darkblue,urlcolor=darkblue,linkcolor=darkblue}
\allowdisplaybreaks

\begin{document}


\title{Controlling unstable chaos: Stabilizing chimera states by feedback}



\author{Jan~Sieber$^{1}$, Oleh~E.~Omel'chenko$^{2,3}$, and Matthias~Wolfrum$^{2}$}
\affiliation{$^1$College of Engineering, Mathematics and Physical Sciences, University of Exeter, North Park Road, Exeter EX4 4QF, United Kingdom}%
\affiliation{$^2$Weierstrass Institute, Mohrenstrasse 39, 10117 Berlin, Germany}%
\affiliation{$^3$Institute of Mathematics, National Academy of Sciences of Ukraine, Tereschenkivska Street 3, 01601 Kyiv, Ukraine}%
\date{\today}

\begin{abstract}
\noindent
We present a control scheme that is able to find and stabilize
an unstable chaotic regime in a system with a large number of interacting particles.
This allows us to track a high dimensional chaotic attractor
through a bifurcation where it loses its attractivity.
Similar to classical delayed feedback control, the scheme is non-invasive,
however, only in an appropriately relaxed sense considering the chaotic regime
as a statistical equilibrium displaying random fluctuations as a finite size effect.
We demonstrate the control scheme for so-called chimera states,
which are coherence-incoherence patterns in coupled oscillator systems.
The control makes chimera states observable close to coherence,
for small numbers of oscillators, and for random initial conditions.
\end{abstract}


\maketitle

\paragraph*{Introduction.}
The classical goal of control is to force a given system
to show robustly a behavior a-priori chosen by the engineer (say, track a desired trajectory).
However, feedback control can also be an analysis tool in nonlinear dynamics:
whenever the feedback input $u(t)$ is zero, i.e the control is \emph{non-invasive},
one can observe natural but dynamically unstable regimes of the uncontrolled nonlinear system
such as equilibria or periodic orbits~\cite{ss2007}.
A famous example is the method of time-delayed feedback control~\cite{p1992},
which provides a non-invasive stabilization of unstable periodic orbits and equilibria~\cite{h2011}.
In general, a control scheme can be useful for nonlinear analysis
if the controlled system converges to an invariant set of the uncontrolled system
without requiring particular a-priori knowledge about the location of the invariant set.
In this context the term ``chaos control'' is used to describe the stabilization
of an unstable periodic orbit that is embedded into a chaotic attractor.
Thus, classical chaos control refers to suppressing chaos~\cite{ogy1990,ss2007}.

In this Letter, we present a control scheme that is able to stabilize
a high-dimensional chaotic regime in a system with a large number of
interacting particles.  Our example is a so-called {\it chimera state},
which is a coherence-incoherence pattern in a system of
coupled oscillators.  We demonstrate that at its point of
disappearance this chaotic attractor turns into a chaotic saddle,
which in our numerical simulation we are able to track as a stable
object by applying the control scheme.  The control scheme is a
classical proportional control that acts globally on a spatially
extended system, as has been used, e.g., for the control of
reaction-diffusion patterns~\cite{ms2006}.  For a chaotic regime,
control is \emph{non-invasive on average} in the following sense: (i) $\langle
u\rangle\to0$ for $t\to\infty$: the time average of the control input
tends to zero over time intervals of increasing length.  (ii) $u\to0$
for $N\to \infty$: the control becomes small for an increasing number
of particles.  The limit $N\to\infty$ has been studied in detail for
chimera states.  Chimera states are stationary solutions of a
well-understood continuum limit system~\cite{oa2008,l2009_,o2013}.
This enables us to compare the chaotic saddle in the finite oscillator
system with the corresponding saddle equilibrium in the continuum
limit system.  However, our control method does not depend on the
knowledge of such a limit and it may be useful in general to
numerically detect a tipping point of a macroscopic state with an
irregular motion on a microscopic level.  On the other hand, we will
show that the proposed control scheme also works for small system
size, where the continuum limit provides only a rough qualitative
description.
\begin{figure}[t]
\includegraphics[width=\columnwidth]{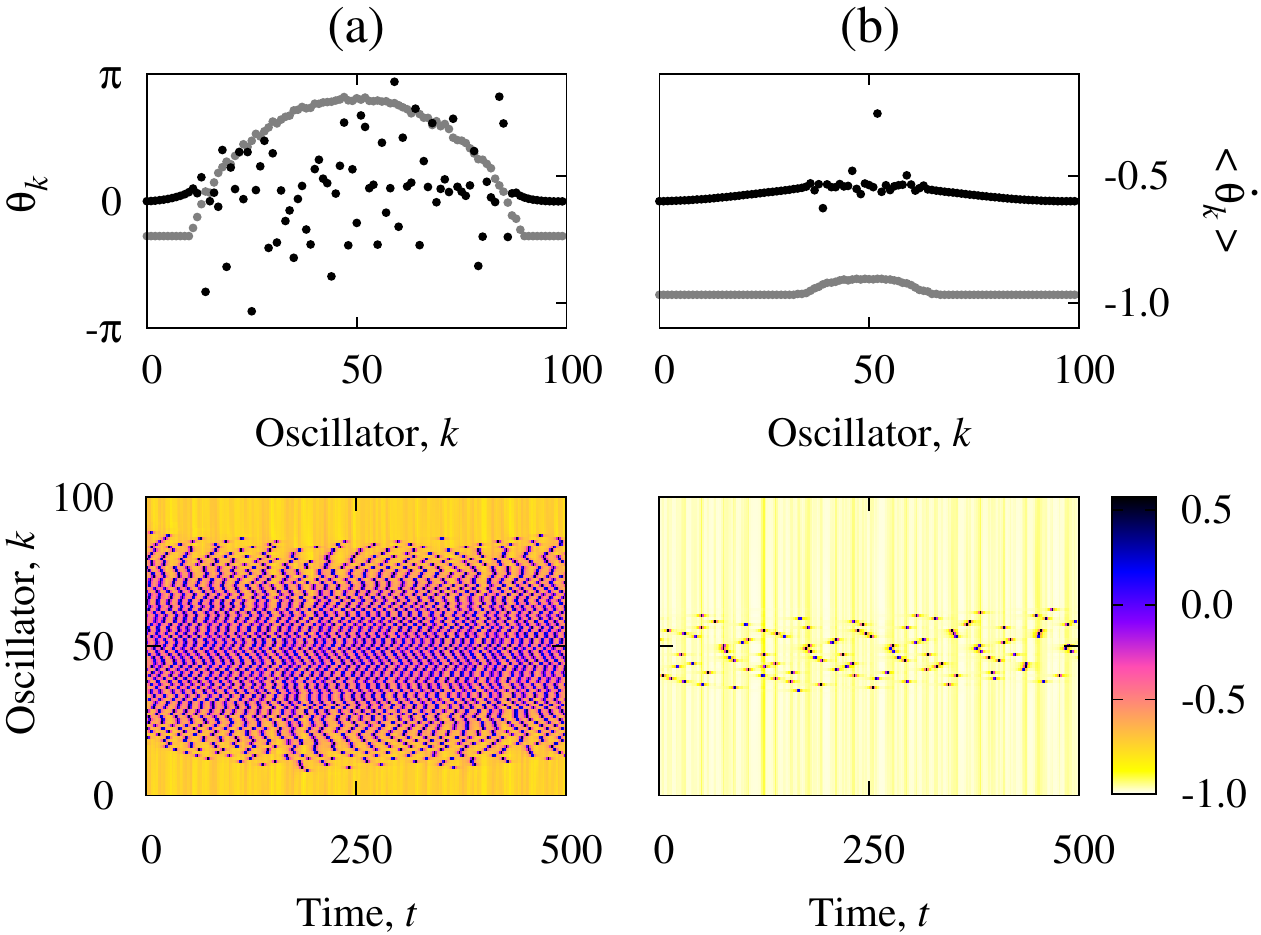}
\caption{(color online) Chimera states far away from complete coherence~(a)
and close to coherence~(b), obtained by numerical simulation
of~\eqref{eq:PhaseOscillators},\,\eqref{eq:gdef} with~$A=0.9$.
Upper panels: Snapshot of phases (black)
and time-averaged phase velocities (gray).
Lower panels: Space-time plots of angular velocities.
We require feedback control~\eqref{eq:control} to observe pattern~(b).}
\label{fig:velocity}
\end{figure}

Applying the control scheme permits us to study the macroscopic state
in regions of the phase and parameter space
that are inaccessible in conventional simulations or experiments.
In the coupled oscillator system this reveals
several interesting properties of the stabilized chimera states.
In the controlled system, we observe a stable branch of chimera states
bifurcating from the completely coherent (synchronized) solution.
This represents a new mechanism for the emergence of a self-organized pattern
from a spatially homogeneous state.
We will show that the dynamical regime of a chimera state
close to complete coherence can be described as a state of {\it self-modulated excitability}.
Moreover, it turns out that also the chimera states on the primarily stable branch
change their stability properties under the influence of the control.
It is known that in the uncontrolled system the chimera states have a dormant instability
that will lead eventually to a sudden collapse of the pattern~\cite{wo2011}.
We will show that this collapse can be successfully suppressed by the control.
Since the chimera's life-span as a chaotic super-transient~\cite{tl2008}
increases exponentially with the system size,
this collapse suppression provides stable chimera states also for very small system size.
In addition to the collapse suppression, the control enlarges the basin of attraction
such that random initial conditions converge almost surely to the chimera state,
which is of particular importance for experimental
realizations~\cite{ttwhs2009,tns2012,hmrhos2012,mtfh2013,nts2013}.

\paragraph*{Chimera states in coupled oscillator systems.}
A chim\-era state is a regime of spatially extended chaos~\cite{woym2011}
that can be observed in large systems of oscillators~\cite{kb2002,as2004} with non-local coupling.
It has the peculiarity that the chaotic motion of incoherently rotating oscillators
is confined to a certain region by a self-organized process of pattern formation
whereas other oscillators oscillate in a phase-locked coherent manner (see Fig.~\ref{fig:velocity}(a)).
The prototypical model of coupled phase oscillators has the form
\begin{equation}
\frac{d \theta_k}{d t} = \omega - \frac{2\pi}{N} \sum\limits_{j=1}^N G_{kj} \sin( \theta_k - \theta_j + \alpha)
\mbox{,}\quad k=1\ldots N
\label{eq:PhaseOscillators}
\end{equation}
where the coupling matrix~$G$ determines the spatial arrangement of the oscillators.
Well-studied cases are rings~\cite{kb2002,as2004,ssa2008,l2009_,woym2011,wo2011,oohs2013},
two-tori~\cite{owyms2012,pa2013} and the plane~\cite{sk2004,mls2010}.
We choose here a ring of oscillators and
\begin{equation}
G_{kj}=G(x_k - x_j) = \frac{1}{2\pi}[1+A\cos(x_k - x_j)]\mbox{,}\label{eq:gdef}
\end{equation}
where $x_k=2k\pi/N-\pi$ is the location of oscillator~$k$ on the ring
and $\theta_k\in[0,2\pi)$ is its phase.
Considering~$x$ as a continuous spatial variable,
one can derive the continuum limit equation
\begin{equation}
\frac{d z}{d t} = \mathrm{i}\omega z + \frac{1}{2}e^{-\mathrm{i}\alpha}\mathcal{G}z - \frac{1}{2}e^{\mathrm{i}\alpha} z^2\mathcal{G}\overline{z}
\label{eq:contlimit}
\end{equation}
for the complex local order parameter~$z(x,t)$, see~\cite{oa2008,l2009_,o2013} for details.
The non-local coupling is here given by the integral convolution
$$
(\mathcal{G}\varphi )(x):=\int_{-\pi}^\pi G(x-y)\varphi (y) \mathrm{d} y.
$$
In this limit a chimera state is represented by a  uniformly rotating solution of the form
\begin{equation}
z(x,t)= a(x)e^{\mathrm{i}\Omega t},
\label{eq:solution}
\end{equation}
where~$\Omega$ is a constant frequency
and~$a(x)$ is a time-independent non-uniform spatial profile
including coherent regions characterized by $|a(x)|=1$
and incoherent regions where  $|a(x)|<1$, see e.g. \cite{o2013}.
\begin{figure}[t]
\includegraphics[width=0.9\columnwidth]{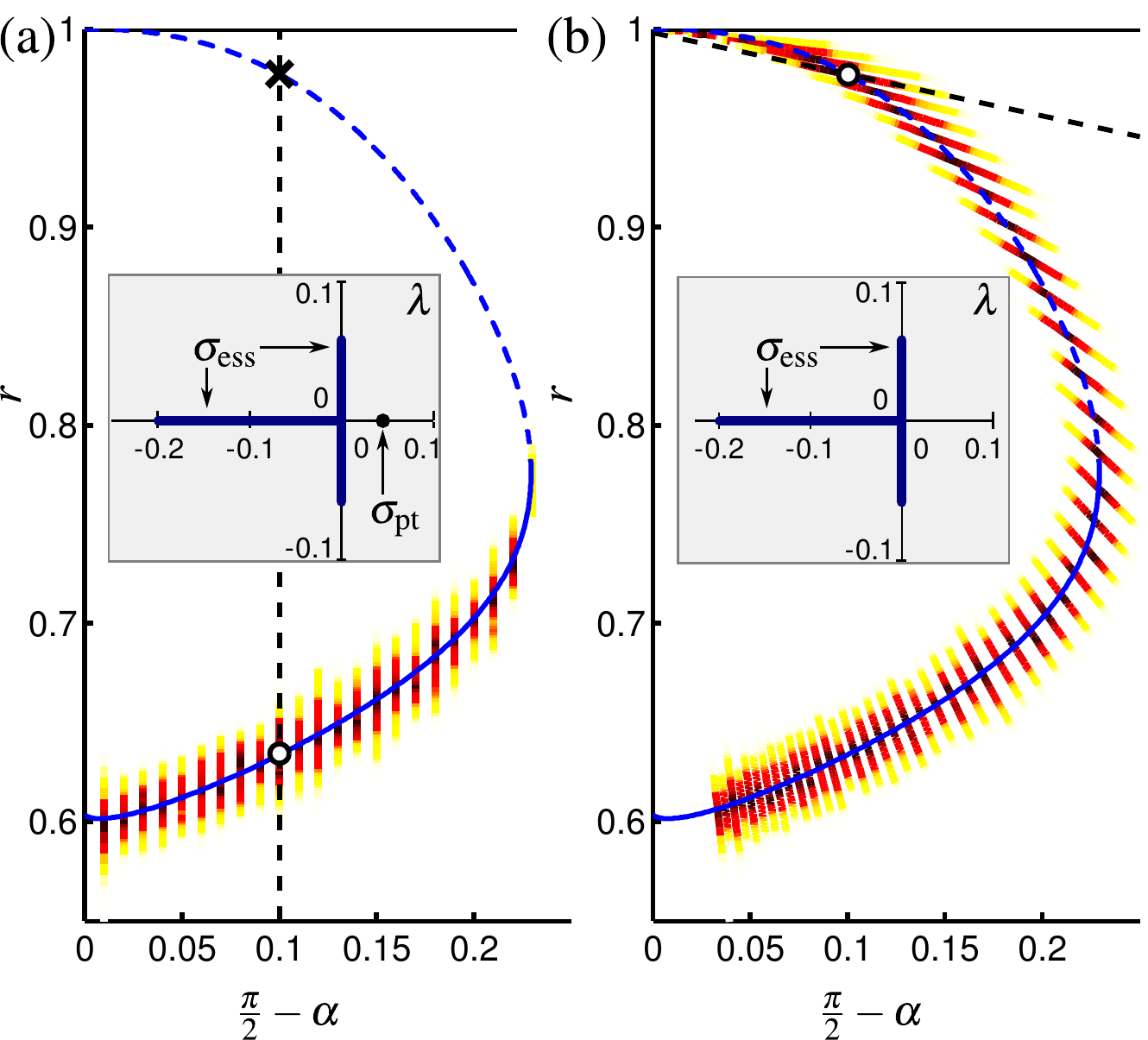}
\caption{(color online) Chimera states projected to the~$(\alpha,r)$
  plane ($N=400$, $A=0.9$).  Panel~(a): Uncontrolled chimeras;
  sequence of simulation runs with stepwise decreasing
  parameter~$\alpha$.  Panel~(b): Controlled chimeras; sequence of
  simulation runs with stepwise increasing control gain~$K$.  Blue
  curve: numerically computed chimera solution of the continuum limit
  (solid: stable; dashed: unstable).  Color/shade patterns:
  observed density in each run (darker$=$higher density, see also
  histograms in Fig.~\ref{fig:invasiveness}).  Highlighted runs along
  dashed lines correspond to the parameter values used in
  Fig.~\ref{fig:velocity} and Fig.~\ref{fig:invasiveness}.
    Insets: spectra of the linearized continuum limit
    \eqref{eq:linearized} for corresponding unstable (a) and
    stabilized (b) chimera state, marked at
    $(\alpha,r)=(\pi/2-0.1,0.98)$.}
\label{fig:biflines}
\end{figure}

A chimera state with finite~$N$ shows temporal and spatial fluctuations
around the corresponding stationary limiting profile.
The color/shade patterns in Fig.~\ref{fig:biflines}(a)
show the stationary densities of the global order parameter
\begin{displaymath}
 r(t) = \frac{1}{N} \left| \sum_{k=1}^N e^{\mathrm{i} \theta_k(t)} \right|
\end{displaymath}
fluctuating around its mean value for a series of chimera trajectories
with varying parameter~$\alpha$.  For the continuum limit
\eqref{eq:contlimit} we obtain a continuous branch of chimera
solutions \eqref{eq:solution} shown as a blue curve in
Fig.~\ref{fig:biflines}, using the continuum version
\begin{equation}
r(t) = \frac{1}{2\pi} \left| \int_{-\pi}^\pi z(x,t) \mathrm{d} x \right|
\label{eq:r_cont}
\end{equation}
for the global order parameter, which is constant for a chimera state~\eqref{eq:solution}.
As Fig.~\ref{fig:biflines}(a) shows, the chimera state disappears
if one decreases the parameter~$\alpha$ beyond~$\pi/2-0.22$.
In the context of the continuum limit~$N\to\infty$ this corresponds
to a classical fold of the solution branch, which continues as an unstable solution
up to the completely coherent state at $(\alpha=\pi/2,r=1)$.

\begin{figure}[t]
\includegraphics[width=\columnwidth]{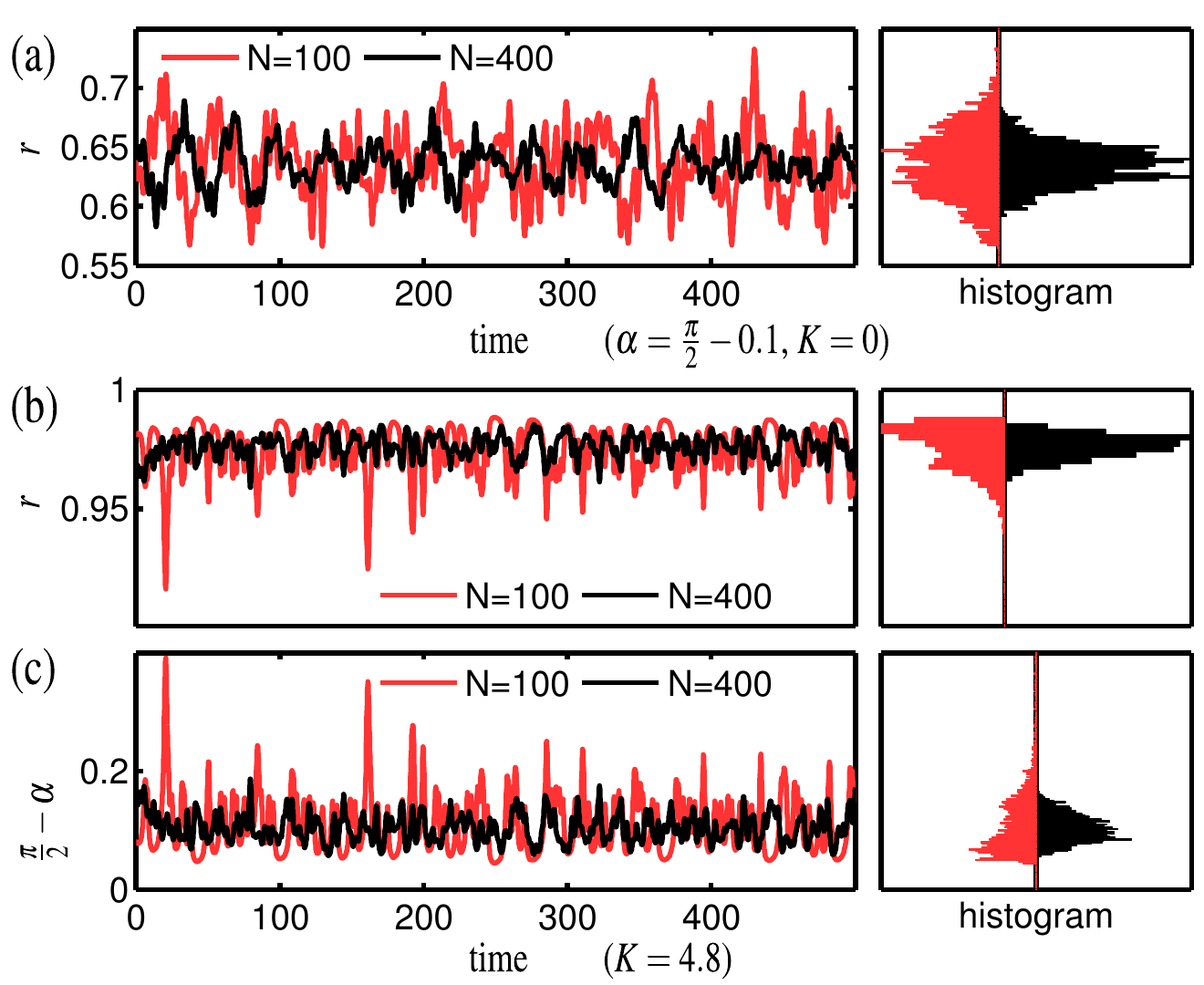}
\caption{(color online)  Time profiles and histograms of global order parameter~$r$
for chimera without control~(a), and~$r$ and~$\alpha$
for chimera with feedback control~(b) and~(c),
for~$N=100$ and~$N=400$ oscillators ($K=4.8$ for~(b,c), $A=0.9$).}
\label{fig:invasiveness}
\end{figure}

\paragraph*{Control scheme.}
In order to study this unstable branch in more detail for moderately sized~$N$
without relying on the continuum limit,
we employ the proportional control scheme
\begin{equation}\label{eq:control}
  \alpha(t)=\alpha_0+K(r(t)-r_0)\mbox{,}
\end{equation}
where the reference point $(\alpha_0,r_0)$ and the control gain~$K$
determine a straight line in the $(\alpha,r)$-plane along which the
controlled system evolves in time (see dashed lines in Fig.~\ref{fig:biflines}).
Setting $K=0$ corresponds to a vertical line, $K\to\infty$ to a horizontal line.
In Fig.~\ref{fig:biflines}(b) we show a sequence of stationary
densities for chimera states in the plane~$\pi/2-\alpha$ vs. global
order parameter~$r$, obtained from numerical simulations
of~(\ref{eq:PhaseOscillators}), now with control~(\ref{eq:control}),
increasing the control gain~$K$ in steps.  The reference point has
been fixed to $(\alpha_0,r_0)= (\pi/2 +0.01,1)$.  In this way, we find
stabilized chimera states along the whole branch of equilibria from
the continuum limit.  Fig.~\ref{fig:invasiveness} shows in more detail
the invasiveness of the control for the runs highlighted in
Figs.~\ref{fig:biflines}(a) and~\ref{fig:biflines}(b) by the dashed
lines.  Whereas for the uncontrolled run the global order
parameter~$r$ fluctuates around its equilibrium value
from the continuum limit (Fig.~\ref{fig:invasiveness}(a)), in the
controlled run both~$r$ and~$\alpha$ fluctuate around their mean
values (Figs.~\ref{fig:invasiveness}(b) and~(c)).  These fluctuations
decrease for an increasing number of oscillators (compare histograms
for~$N=100$ and~$N=400$ in Fig.~\ref{fig:invasiveness}).  Since for a
finite~$N$ system the invasiveness of the control is given by the
fluctuations of these global quantities, it is non-invasive on average
satisfying conditions (i)--(ii) stated above.

Note that chimera states in a system with a nonlinear state-dependent
phase-lag parameter have been investigated already in~\cite{bpr2010}.
However, the feedback in \cite{bpr2010} depends on the local order
parameter such that it cannot be interpreted as a global non-invasive
control of the original system in the sense of \cite{ms2006}.  Proportional
control~\eqref{eq:control} is only one option to achieve non-invasive
control on average for a chaotic saddle in the relaxed sense of conditions
(i)--(ii).  Alternatives are any non-invasive methods for
stabilization of unknown equilibria.  For example, a PI
(proportional-integral) control was used in~\cite{tmk2010} to explore
the saddle-type branch of a partially synchronized regime in a
small-world network in the continuum limit.
Time-delayed feedback or wash-out filters~\cite{AWC94}
are suitable near instabilities other than folds of the continuum-limit equilibrium;
for instance in~\cite{rp2004}, time-delayed feedback has been used
to suppress or enhance synchronization in a system of globally coupled oscillators.

\paragraph*{Spectral stability analysis.}
In the continuum limit~(\ref{eq:contlimit}), the control~(\ref{eq:control}), (\ref{eq:r_cont})
acts in an exactly non-invasive manner and the stabilization can be shown as follows.
For a solution~(\ref{eq:solution}) with $\alpha = \alpha_\ast$ and $r=r_\ast$, we insert
$$
z(x,t) = ( a(x) + v(x,t) ) e^{\mathrm{i} \Omega t}
$$
into Eq.~(\ref{eq:contlimit}) with control~(\ref{eq:control}), (\ref{eq:r_cont})
and linearize the result with respect to the small perturbation~$v$.
As a result, we obtain the linear equation (c.f.~\cite{woym2011})
\begin{equation}
\frac{d v}{d t} = \mathcal{L} v := \eta(x) v(x,t) + \mathcal{K} v + \mathcal{C} v,
\label{eq:linearized}
\end{equation}
containing the multiplication operator
\begin{equation}
\eta(x) := \mathrm{i} ( \omega - \Omega ) - e^{\mathrm{i} \alpha_\ast} a(x) \mathcal{G}\overline{a}
\end{equation}
and the compact integral operators
\begin{eqnarray*}
(\mathcal{K} v)(x) &:=& \frac{1}{2} e^{-\mathrm{i} \alpha_\ast} \mathcal{G} v - 
\frac{1}{2} e^{\mathrm{i} \alpha_\ast} a^2(x) \mathcal{G} \overline{v},\\[2mm]
(\mathcal{C} v)(x) &:=& \frac{\mathrm{i} K a(x)\eta(x)}{4\pi^2 r_\ast}
\mathrm{Re}\left( \int_{-\pi}^\pi \overline{a} (y) \mathrm{d}y\: \int_{-\pi}^\pi v(y) \mathrm{d}y\right),
\end{eqnarray*}
where~$\mathcal{C} v$  accounts for the action of the control.
Spectral theory for this type of operators (see~\cite{o2013} for details) implies
that the spectrum~$\sigma(\mathcal{L})$ consists of  two qualitatively different parts:
(i) essential spectrum
$$
\sigma_{\mathrm{ess}}(\mathcal{L}) = \{ \eta(x)\::\: -\pi\le x\le\pi \}\cup\{ c.c. \},
$$
which for partially coherent states is known to have a neutral part \cite{ms2007};
(ii) point spectrum~$\sigma_{\mathrm{pt}}(\mathcal{L})$
consisting of all isolated eigenvalues of the operator~$\mathcal{L}$.
For the chimera states shown in Fig.~\ref{fig:biflines},
the point spectrum contains at most one real eigenvalue,
which determines their stability.
This eigenvalue can be found by  inserting $v = v_0(x) e^{\lambda t}$ into Eq.~(\ref{eq:linearized}),
\begin{equation}
v_0(x) = ( \lambda - \eta(x) )^{-1} ( \mathcal{K} v_0(x) + \mathcal{C} v_0(x) ).
\label{eq:EM}
\end{equation}
Applying now the integral operator~$\mathcal{K} +  \mathcal{C}$ to both sides of Eq.~(\ref{eq:EM})
we arrive at a spectral problem for $w := (\mathcal{K} +  \mathcal{C}) v_0$
\begin{equation}
w = (\mathcal{K} +  \mathcal{C}) \left( ( \lambda - \eta(x) )^{-1} w \right).
\label{eq:EM_Transformed}
\end{equation}
As pointed out in~\cite{as2004}, the operators~${\cal G}$ and~${\cal  K}$
have finite rank for the coupling function~\eqref{eq:gdef}
(the control term ${\cal C}$ has always rank one).
Therefore, expanding~$w$ as a Fourier series and projecting Eq.~(\ref{eq:EM_Transformed})
onto the first three modes $f_1(x) = 1$, $f_2(x) = \cos x$, $f_3(x) = \sin x$,
we obtain a closed linear system
\begin{equation}
\hat{w}_k = \frac{1}{\pi} \int_{-\pi}^\pi f_k(x)  (\mathcal{K} +  \mathcal{C}) \left( ( \lambda - \eta(x) )^{-1} w \right) \mathrm{d}x
\label{eq:EM_Fourier}
\end{equation}
for the unknown Fourier coefficients~$\hat{w}_0$, $\hat{w}_1$
and~$\hat{w}_2$.  Its determinant gives an equation that is nonlinear
for the real eigenvalues~$\lambda$ and linear in the gain $K$.
The insets in Fig.~\ref{fig:biflines} show the spectra calculated in this way,
indicating the unstable eigenvalue in panel~(a), which disappears due to the control~(b).

\begin{figure}[t]
\includegraphics[width=\columnwidth]{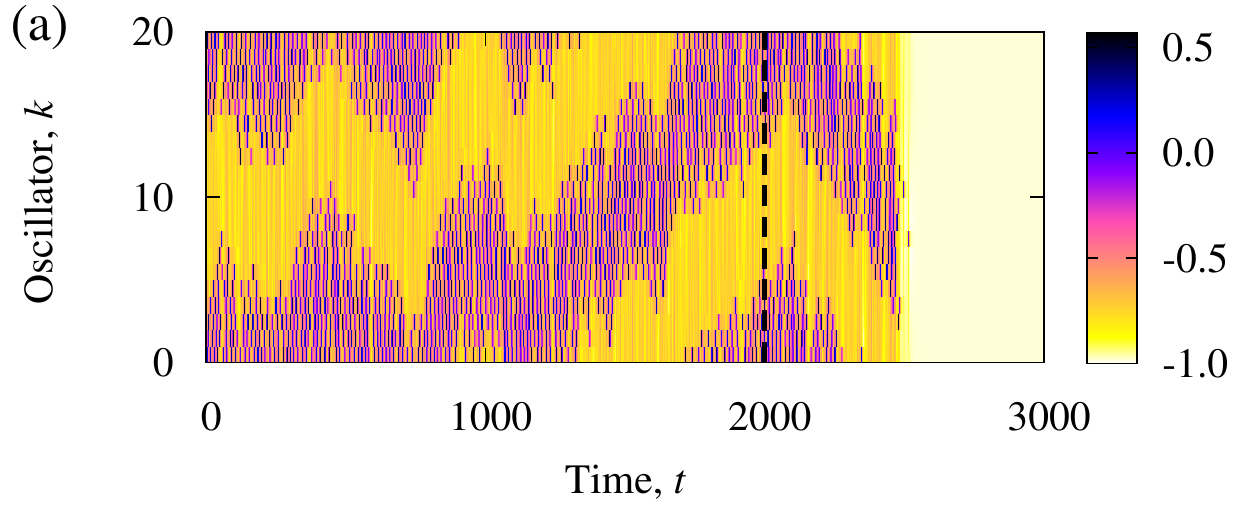}
\includegraphics[width=\columnwidth]{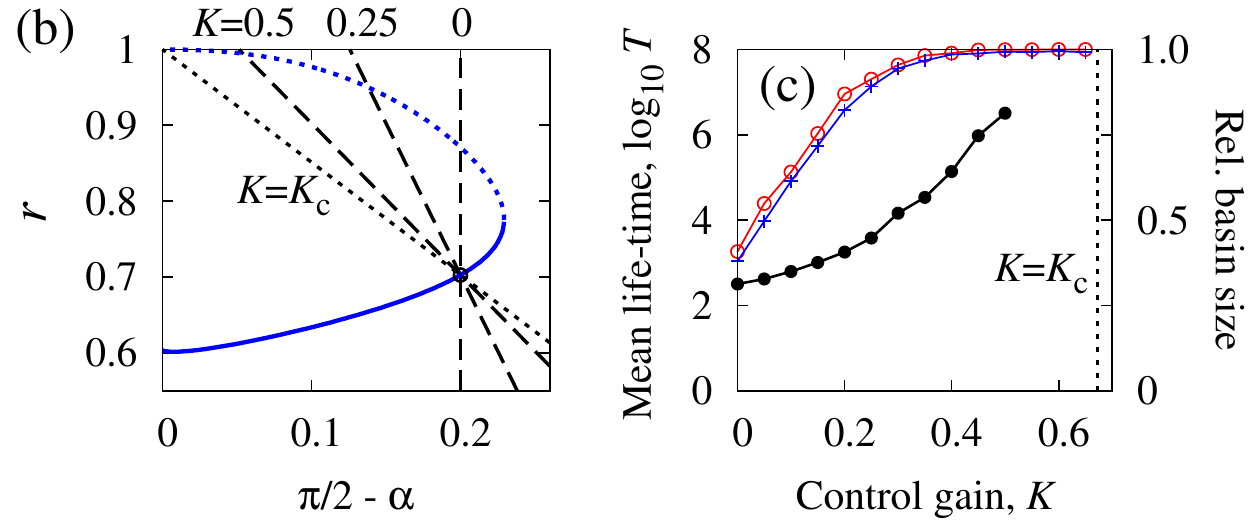}
\caption{(color online) Influence of the control on a stable chimera state.
Panel~(a): Space-time plot of angular velocities;
switching off the control with~$K=1$ at $t=2000$
permits the subsequently observed collapse for~$N=20$.
Panel~(b): Controlling the same chimera state with increasing values of the control gain~$K$.
Panel~(c): mean life-time before collapse  for $N=20$ (dots);
fraction of random initial conditions attracted
by the chimera state for~$N=20$ (circles) and~$N=100$ (crosses)
}
 \label{fig:collapse}
 \end{figure}

\paragraph*{Suppression of collapse and enlarged basins.}
We study now the influence of the control scheme
on the classical chimera states far from complete coherence,
which are already stable without the control (solid blue curve in Fig.\ref{fig:biflines}(a)).
As described in~\cite{wo2011}, the classical chimera states
from time to time show a sudden transition to the stable completely coherent state
and have to be considered as weakly chaotic type-II supertransients~\cite{tl2008}.
The life-time before collapse increases exponentially with the system size
which implies that chimera states disappear quickly for~$N\approx 20$ (cf. Fig.~\ref{fig:collapse}(a)),
whereas they typically appear as stable objects for any observable time-span if~$N>100$.
The collapse process can be understood as follows.
Driven by finite size fluctuations, the trajectory can tunnel through the barrier
represented by the chimera on the unstable branch
and eventually reach the stable coherent state.
Applying the control, this scenario changes drastically:
Increasing the control gain~$K$, the mean life-time before collapse
increases by several orders of magnitude and, at the same time,
the basin of attraction of the chimera state grows correspondingly.
Fig.~\ref{fig:collapse}(c) shows the average observed life-times for increasing values of~$K$.
In our simulations over~$10^7$ time units, which we performed for each~$K$,
the number of observed collapses decreased successively until for $K>0.5$,
we did not observe a single collapse event during this time span.
Finally, for $K\geq K_{\mathrm c}\approx 0.67$ the chaotic saddle acting as a barrier
disappears and the completely coherent state becomes unstable,
which ultimately prevents a collapse to this state.
Accordingly, all random initial data converged to the chimera state.
Note that we have chosen the reference point on the chimera branch, see Fig.~\ref{fig:collapse}(b),
such that the given chimera state exists for all values of the control gain~$K$.
Hence, with feedback control stable chimera states can be observed
for considerably smaller values of~$N$, and arbitrary initial conditions,
which is of particular importance for experimental realizations.
\paragraph*{Self-modulated excitability close to coherence.}
Up to now, stable chimera states have been observed
only far from the completely coherent solution, except for the results in~\cite{omt2008}
where the onset of incoherence has been triggered by an inhomogeneous stimulation profile.
In the controlled system~\eqref{eq:PhaseOscillators}, \eqref{eq:control}
there is a stable branch of chimera states bifurcating from complete coherence.
This is another example of a pattern forming bifurcation mechanism in a homogeneous system
with a diffusion like coupling that should in principle stabilize homogeneity.
The chimera states close to complete coherence display particular properties
distinguishing them from classical chimera states.
Fig.~\ref{fig:velocity}(b) shows that the onset of incoherence manifests itself
as the emergence of isolated excitation bursts caused by phase slips of single or few oscillators,
which appear irregular in space and time but are confined by a process of self-localization  to a certain region.
Indeed, close to the bifurcation point the dynamics of each single oscillator
is close to a saddle-node-on-limit-cycle bifurcation.
Hence, the emergence of a chimera state can be understood
as a transition from quiescent to oscillatory behavior,
which happens in a self-localized excitation region within a discrete excitable medium.
At the same time, the isolated phase slipping events are not well described
by the average quantities from the continuum limit,
which are continuous in space and constant in time.

\paragraph*{Conclusion.}
We demonstrate that a feedback control that is non-invasive in our
relaxed sense is useful for exploring complex dynamical regimes in
large coupled systems.  In particular, it can be used to classify the
disappearance of a chaotic attractor as a transition to a chaotic
saddle, which is the classical scenario for so-called \emph{tipping},
e.g., in climate \cite{lhkhlr2008}, without relying on a closed-form
continuum limit.  Specific to partial coherence, feedback control is
feasible and useful in existing experimental setups of coupled
oscillators~\cite{mtfh2013,tns2012,hmrhos2012,nts2013} as the coupling
in these experiments is computer controlled or through a mechanical
spring. Feedback control makes it possible to study the phenomenon of
partial coherence for much smaller~$N$, close to complete coherence,
and without specially prepared initial conditions.

%

\end{document}